\documentclass[draft]{amsart} 
\usepackage{amssymb,latexsym,amsfonts,verbatim,amscd} 

\numberwithin{equation}{section}

\theoremstyle{plain}

\newtheorem{theorem}[equation]{Theorem}   
\newtheorem{lemma}[equation]{Lemma} 
\newtheorem{proposition}[equation]{Proposition} 
\newtheorem{corollary}[equation]{Corollary}

\theoremstyle{definition}

\newtheorem{definition}[equation]{Definition} 
\newtheorem{remark}[equation]{Remark} 
 
\newtheorem{example}[equation]{Example}

\DeclareMathOperator{\Hom}{Hom}

\DeclareMathOperator{\HH}{H}

\DeclareMathOperator{\coker}{Coker}

\DeclareMathOperator{\rank}{rank} 
\DeclareMathOperator{\im}{Im}

\DeclareMathOperator{\lcm}{lcm} 
\DeclareMathOperator{\Ker}{Ker} 
\DeclareMathOperator{\sgn}{sgn}

\def\mapright#1{\smash{\mathop{\longrightarrow}\limits^{#1}}}

\renewcommand\a{\alpha} 
\renewcommand\b{\beta}

\renewcommand{\:}{\! :}

\newcommand{\lra}{\longrightarrow}

\newcommand{\smsm}{\smallsetminus}

\newcommand{\e}{\epsilon} 
\newcommand{\g}{\gamma} 
\newcommand{\U}{\upsilon}

\begin{document}

\title[Free resolutions for multigraded modules]
{Free resolutions for multigraded modules: \\ 
       a generalization of Taylor's construction} 
\author[H. Charalambous]{Hara Charalambous} 
\address{Department of Mathematics\\ 
         University at Albany, SUNY\\ 
         Albany, NY 12222} 
\email{hara@math.albany.edu} 
\email{tchernev@math.albany.edu} 
\author[A. Tchernev]{Alexandre Tchernev} 
\date{\today}

\begin{abstract} 
Let $Q=\Bbbk[x_1,\ldots, x_n]$ be a polynomial ring over a field 
$\Bbbk$ with the standard $\Bbb N^n$-grading. 
Let $\phi$ be a morphism of finite free $\mathbb N^n$-graded 
$Q$-modules. We translate to this setting several notions and 
constructions that appear originally in the context of monomial 
ideals. First, using a modification of the Buchsbaum-Rim complex, 
we construct a canonical complex $T_\bullet(\phi)$ of  finite free 
$\Bbb N^n$-graded $Q$-modules that generalizes Taylor's resolution.  
This complex provides a free resolution 
for the cokernel $M$ of $\phi$ when $\phi$ satisfies certain rank criteria. 
We also introduce the Scarf complex of $\phi$, and a notion of 
``generic'' morphism. Our main result is that the Scarf complex of 
$\phi$ is a minimal free resolution of $M$ when $\phi$ is minimal and generic.   
Finally, we introduce the LCM-lattice for $\phi$ and establish its 
significance in determining the minimal resolution of $M$.
\end{abstract}

\maketitle


\tableofcontents 

\section{Introduction}\label{S: introd}

There is a plethora of significant papers examining  
free and minimal resolutions 
of monomial ideals. In contrast relatively little is known for 
$\mathbb N^n$-graded (multigraded) modules. 
Historically the prototype of a free resolution   
for monomial ideals is the 
Taylor resolution, \cite{Ta60}. More recently \cite{BaPeSt98} 
give the minimal free 
resolution of generic monomial ideals, based on the idea of the 
Scarf complex. Soon after  that \cite{GaPeWe99} discuss the
significance of the LCM-lattice in determining 
the minimal free resolution of an ideal.  For multigraded modules, 
\cite{ChDe01} discuss the second syzygies, while \cite{Ya99} 
and \cite{Mi00} among others generalize 
results concerning 
homological invariants of monomial ideals to multigraded modules. 

In this paper we translate to the setting of morphisms of finite 
free multigraded modules several notions and 
constructions that appear originally in the context of monomial 
ideals. For a morphism $\phi$ of finite free multigraded modules, 
by using the formalism of what we 
decided to call \emph{Buchsbaum-Rim-Taylor systems}, we 
construct canonical complexes $T_\bullet(\phi)$ and $S_\bullet(\phi)$ 
of  finite free multigraded modules. 
The \emph{Taylor complex} $T_\bullet(\phi)$ generalizes Taylor's 
resolution of a monomial ideal,  and provides a free resolution 
for the cokernel $M$ of $\phi$ when certain rank 
criteria are satisfied for $\phi$. 
Just as the underlying linear algebra of the Taylor 
resolution is based on the Koszul complex, the 
linear algebra structure of the complex $T_\bullet(\phi)$ is based 
on the Buchsbaum-Rim complex \cite{BuRi63}.  
The \emph{Scarf complex} $S_\bullet(\phi)$ is a minimal subcomplex 
of the Taylor complex, and appears to be the appropriate 
generalization of the Scarf complex of a monomial ideal. 

We also introduce a  
notion of ``generic'' morphism. It is based on the combinatorial 
notion of generic monomial ideal of \cite{BaPeSt98}, together 
with the requirement that the morphism be sufficiently generic also 
from linear algebra point of view. 
Our main result, Theorem \ref{T:Scarfres}, is that the Scarf complex of 
$\phi$ is a minimal free resolution of $M$ when $\phi$ is a minimal 
presentation of $M$ and is a generic morphism.   
Finally, we introduce the LCM-lattice for a morphism 
$\phi$, and establish its 
significance in determining the minimal resolution of $M$.

The paper is organized as follows. 
In Section \ref{S: prelim} we set our notation.  
For a morphism $\phi$ we introduce the 
coefficient matrix: a matrix with entries from 
$\Bbbk$ which determines the linear algebra structure of the 
resolution. We define certain submatrices of the coefficient 
matrix, and the kernel of a dual map which will turn out essential 
for the proof of the main theorem of this paper, Theorem 
\ref{T:Scarfres}.

In Section \ref{S:BRcomplex} we recall the construction of the 
Buchsbaum-Rim complex , presenting  it in a form 
that guarantees that the differentials of the complex stay invariant 
under change of basis.  This is important because, as it will be clear 
in Section \ref{S: minimalres}, to compute the minimal resolution one has 
to consider a change of basis depending on the kernel of the dual map
of  Section \ref{S: prelim}.

In Section \ref{S: BRTtheory} we introduce the notion of a BRT 
(Buchsbaum-Rim-Taylor) system, and use it to define a  
BRT complex for $\phi$.  The BRT complex for the full  BRT system is 
exact when the ranks of the submatrices 
of the coefficient matrix are high enough. This full BRT complex is the 
Taylor complex. 
 
In Section  \ref{S: minimalres} we introduce the LCM-lattice, the 
Scarf simplicial complex, and the notion of generic morphism. 
We describe the
BRT system whose BRT complex is the Scarf complex and formulate  
our main result, Theorem \ref{T:Scarfres}. 

Section \ref{S:Proof} is devoted to the proof of Theorem 
\ref{T:Scarfres}. 
 
Finally, 
in Section \ref{S: LCM} we show how to extend to 
a certain class of morphisms (the morphisms of \emph{uniform rank}) 
the arguments in \cite{GaPeWe99} Theorem 3.3.  This allows us 
to exhibit the role 
that the LCM-lattice plays in determining the structure of the minimal 
resolutions of these morphisms.

\section{Preliminaries}\label{S: prelim}  

Throughout this paper $\Bbbk$ is a field, and 
$Q = \Bbbk[x_1,\ldots,x_n]$ is  the 
polynomial ring in $n$ variables over $\Bbbk$. 
Let $\a=(a_1,\dots, a_n)$ be an element of $\mathbb N^n$. 
The \emph{support} of $\a$ is the set $supp(\a)=\{i\mid a_i\ne 0\}$. 
We write $x^{\a}$ for the monomial $x_1^{a_1}\dots x_n^{a_n}$, and 
we set the degree of $x^{\a}$ to be $|x^{\a}|=\a$. 
This makes  $Q$ into a $\mathbb N^n$-graded (or {\it multigraded}) 
algebra. We consider the partial order on  $\mathbb N^n$ 
given  by 
\[
\a=(a_1,\dots, a_n) \preceq \b=(b_1,\dots, b_n) 
  \quad \iff \quad a_i\le b_i \text{ for } i=1,\dots,n. 
\]
In addition, we define the {\it join} or {\it lcm} of $\a$ and $\b$ by 
\[
\lcm(\a, \b)= \a \vee \b =\bigl(\max(a_1,b_1),\dots,\max(a_n,b_n)\bigr). 
\]
The tensor product of multigraded $\Bbbk$-vector spaces is 
multigraded with $|x\otimes y|=|x| + |y|$. Unadorned tensor 
products are over $\Bbbk$. 
 
Let $E$ and $G$ be finite free multigraded $Q$-modules of ranks 
$e$ and $g$ respectively, and let $\phi\: E \lra G$ be a 
multigraded morphism. We fix homogeneous bases 
$\e_1,\dots, \e_e$ of $E$ and $\g_1,\dots, \g_g$ of $G$,  
and we let $\Phi=(f_{ij})$ be the matrix of $\phi$ in these bases.  
Thus the $j$th column $\Phi_j$ 
of $\Phi$ gives the image of $\e_j$ in $G$. 
We say that the degree of the column $\Phi_j$ is $|\e_j|$. 

To every such map $\phi$ we associate a map $s\: U\lra W$ of 
vector spaces $U$ and $W$, where we let $U$ be 
the $\Bbbk$-vector space with basis $\e_1,\dots,\e_e$ and 
let $W$ be the $\Bbbk$-vector space with basis $\g_1,\dots,\g_g$   
The matrix $C$ of the map $s$ in terms of the given bases is the 
\emph{coefficient matrix} of $\Phi$\ : 
each entry $f_{ij}$ of $\Phi$ is of the form 
$f_{ij}=c_{ij}x^{\a_{ij}}$ where $c_{ij}\in\Bbbk$, and $C=(c_{ij})$. 
Clearly  $\rank C=\rank \Phi= \rank \phi$.

Let   $\a=(a_1,\dots,a_n)$ be a multidegree. In Section 
\ref{S: minimalres} we will make use of
the maps $s_{\a}$ and the vector spaces $V_{\a}$ and $K_{\a}$ 
which we define below. 
First, let $U_{\a}$  be the vector subspace of $U$ with 
basis those basis vectors of $U$ whose multidegree in $E$ is 
at most $\a$. Next, let 
$\Phi_{\a}$  be the submatrix of $\Phi$ on columns of degree 
$\preceq \a$, let $C_{\a}$ be the coefficient matrix of 
$\Phi_{\a}$, and let $s_{\a}$  be the restriction of the map 
$s$ to $U_{\a}$. The matrix  of $s_{\a}$ is $C_{\a}$, 
and we write $V_{\a}$ for the image of $U_{\a}$ under $s$. 
Note that $U_{\a}\subseteq U_{\b}$ and 
$V_{\a}\subseteq V_{\b}$ when $\a\preceq\b$.  It is clear 
that the definitions of $U_{\a}$, $V_{\a}$, and $s_{\a}$ 
are independent of the choice of the homogeneous basis of $E$.  

Let $I$ be a subset of $\{1,\dots,e\}$. We write $U_I$ for the 
vector subspace of $U$ with basis $\{\e_i \mid i\in I\}$. We 
call $V_I$ the image of $U_I$ under $s$, and denote by $s_I$ 
the restriction of $s$ to $U_I$. 

For a $\Bbbk$-vector space $Z$ we set 
$Z^*=\Hom_\Bbbk(Z,\Bbbk)$. If $V$ is the image of $s$ then 
the inclusions $V_{\a}\lra V$ and $V_I\lra V$ 
induce surjections $V^*\lra V_{\a}^*$ and $V^*\lra V^*_I$, 
and we write $K_{\a}$ and $K_I$ for the 
corresponding kernels. Note that if $v\in K_\a$ 
and  $u\in U_\a$ then $v\bigl(s(u)\bigr)=0$. 

Finally, if $\mathbf M$ is a complex of vector spaces with 
$i$th differential $\partial_i\: M_i\lra M_{i-1}$, then its  
shift $\mathbf M[k]$ is the complex with $M[k]_i=M_{i+k}$ 
and differential given by $\partial[k]_i=\partial_{k+i}$. 
We say that the complex $\mathbf M$ is exact if 
$\HH_i(\mathbf M)=0$ for $i\ne 0$. We say that $\mathbf M$ 
is split exact if it is exact, and $\HH_0(\mathbf M)=0$ 
as well.

\section{The Buchsbaum-Rim complex}\label{S:BRcomplex}

Let $s\: U \lra W$ be a $\Bbbk$-vector space map. 
For the convenience of the reader, and to establish notation, 
we recall in this section the Buchsbaum-Rim complex 
of the map $s$, cf. \cite{BuRi63}, or 
\cite{Ei97} Section A2.6.1. 

Let $V$ be a subspace of rank $r$ of $W$ such that $V$ contains 
the image of $s$. 
We recall that the $i$th divided power $D_iV^*$ is the dual 
$(S_iV)^*$ of the $i$th symmetric power $S_iV$, and refer to 
\cite{Ei97} for the properties of these functors.   
For any integers $m,k\ge 0$ let $A^{m,k}_\bullet(s, V)$ be the complex  
\[
\begin{CD} 
                0 \lra A^{m,k}_{e-k}   @> \sigma^{m,k}_{e-k}   >> 
                       A^{m,k}_{e-k-1} @> \sigma^{m,k}_{e-k-1} >> 
                           \cdots      @> \sigma^{m,k}_2       >> 
                       A^{m,k}_1       @> \sigma^{m,k}_1       >>  
                       A^{m,k}_0\lra 0 
\end{CD} 
\]
where 
\[ 
A^{m,k}_i = D_{m+i} V^* \otimes \wedge^{k+i} U.   
\] 
The differential $\sigma^{m,k}_i$ of $A^{m,k}_\bullet$ is 
defined as the composition 
\[ 
\begin{CD} 
D_{m+i} V^*\otimes\wedge^{k+i} U                            \\ 
@VV {\delta\otimes\delta} V                                 \\  
D_{m+i-1} V^* \otimes V^* \otimes U\otimes\wedge^{k+i-1}U   \\ 
@VV {1\otimes \mu\circ(1\otimes s)} \otimes 1V              \\ 
D_{m+i-1} V^* \otimes \wedge^{k+i-1} U,   
\end{CD} 
\] 
where $\delta$ is the diagonal map, and 
$\mu\: V^*\otimes V\lra \Bbbk$ is the canonical pairing. 

The complexes $A^{0,k}_\bullet(s, V)$ have been extensively studied.  
In the sequel we will use the following well known  property. 

\begin{proposition}\label{P:koszulexact} 
The complex $A^{0,k}_\bullet(s, V)$ is exact if and only if either 
$k\ge e$, or $V=\im(s)$. 
The complex $A^{m,0}_\bullet(s,V)$ is split exact if $m>0$ and $V=\im(s)$. 
\end{proposition}   

\begin{proof} 
The result is clear for $k\ge e$. Assume $0\le k < e$. 
Fixing a basis of $U$ provides an isomorphism 
$\wedge^t U \cong \wedge^{e-t} U^*$ and this 
identifies the complexes $A^{0,k}_\bullet(s, V)$ and 
$A^{m,0}_\bullet(s,V)$ 
with the complexes $C^k$ and $C^{m+e}$ studied by 
Lebelt \cite{Le73}. In particular, 
the proposition is an immediate consequence of 
\cite{Le73}, Corollary 1 to Theorem 5 and Corollary to Theorem 13; 
or see \cite{Tc96}, Theorem 4.1. 
\end{proof}

When $V$ is equal to the image of $s$ we write $A^k_\bullet(s)$ 
instead of $A^{0,k}_\bullet(s, V)$. In that case  
the induced by $s$ map $U\lra V$ is surjective, its dual  
$V^*\lra U^*$ is an inclusion, and we use it to identify $V^*$ 
as a subspace of $U^*$. It is now an elementary exercise in 
multilinear algebra to show that 
$\coker(\sigma^k_1)$ is isomorphic to $\wedge^{e-k}(U^*/V^*)$ 
when $k\le e$, and is $0$ otherwise. In particular, we have  

\begin{corollary}\label{T:splitexact} 
If $k < r$ the complex $A^k_\bullet(s)$ is split exact. 
If $s$ is injective, then the complex $A^k_\bullet(s)$ is 
split exact for $k\ne e$. 
\qed
\end{corollary}

In order to define the Buchsbaum-Rim complex, we splice  
together the complex $A^{0,r+1}_\bullet(s,V)\otimes\wedge^r V^*$ and 
the complex $U\mapright{s}W$, to obtain a diagram 
\[ 
\begin{CD} 
B_\bullet(s, V)= \quad\quad 
A^{0,r+1}_\bullet(s,V)\otimes\wedge^r V^* @>{s_2}>>  
                                      U   @> s >>  W. 
\end{CD} 
\]
Thus the complex $B_\bullet(s, V)$ has the form 
\[ 
\begin{CD} 
B_\bullet(s, V)= 0 \lra  
       B_{e-r+1}@> s_{e-r+1} >> B_{e-r} @> s_{e-r} >>  
          \dots @> s_2       >> B_1     @> s       >> B_0 \lra 0  
\end{CD}  
\] 
where $B_0=W$ and $B_1=U$, while $B_i=A^{r+1}_{i-2}\otimes\wedge^r V^*$ 
for $i\ge 2$. The splice map 
\[ 
s_2\: \wedge^{r+1}U \otimes\wedge^r V^*\lra U   
\] 
is defined as the composition 
\[ 
\begin{CD} 
\wedge^{r+1}U \otimes\wedge^r V^* 
@>{\bigl((1\otimes\wedge^r s)\circ\delta\bigr)\otimes 1}>>  
U\otimes\wedge^r V\otimes\wedge^r V^* @>{1\otimes \mu}>>  
U\otimes \Bbbk = U,  
\end{CD} 
\] 
and for $i\ge 3$ we have $s_i=\sigma^{r+1}_{i-2}\otimes 1$.  
We note that the factor $\wedge^r V^*$ is just a copy of $\Bbbk$, 
but is needed in order to make the differentials of $B_\bullet$  
invariant under change of basis in $V$.  We will have to consider 
such a change of basis in Section \ref{S: minimalres}. 

Now that we have an invariant description of the differentials 
$s_i$, we can describe them in terms of basis elements 
as follows. Let $\e_1,\dots,\e_e$ be a basis of $U$, let 
$\g_1,\dots,\g_g$ be a basis of $W$,  let 
$\U_1,\dots,\U_r$ be a basis of $V^*$, and consider the dual basis 
on $V$. Let 
$C=(c_{ij})$ be the matrix of the induced by $s$ 
map $U \lra V$ with
respect to the above bases. Let $\b=(b_1,\dots,b_r)$ 
be a sequence of integers with $b_1+\dots+b_r=p$. 
The elements of the form 
$\U^{(\b)}=\U^{(b_1)}_1\dots \U^{(b_r)}_r$ 
where all the $b_i$s are nonnegative are a basis for $D_p V^*$, 
and the element 
$\U_{[r]}=\U_1\wedge\dots\wedge\U_r$ forms a basis for 
$\wedge^r V^*$.  Also, we set $\U^{(\b)}=0$ if $b_i<0$ for some 
$i$, and we write $\b_i$ for the sequence 
$(b_1,\dots,b_{i-1},b_i-1,b_{i+1},\dots,b_k)$. 

Similarly, the elements 
$\e_I=\e_{i_1}\wedge\dots\wedge\e_{i_q}$  $\bigl($where the 
subset $I=\{i_1,\dots,i_q\}$ with $i_1<\dots<i_q$ ranges over all 
$q$-element subsets of $\{1,\dots,t\}\bigr)$ form a basis of 
$\wedge^q U$. Then we have for $i\ge 3$ 
\[ 
s_i(\U^{(\b)}\otimes\e_I\otimes\U_{[r]}) = 
\sum_{j=1}^r \ \sum_{l\in I} \sgn(l, I\smsm l) \  
                           \U_j\bigl(s(\e_l)\bigr)\ 
                           \U^{(\b_j)}\otimes\e_{I\smsm l} 
                                      \otimes\U_{[r]}, 
\] 
while 
\[ 
s_2\bigl(\e_J\otimes\U_{[r]}\bigr)= \sum_{l\in J} \sgn(l,J\smsm l)\ 
                                       \det(C_{J\smsm l}) \ \e_l, 
\] 
where $|J|=r+1$.  Note that 
if $V=W$ and $\U_1,\dots,\U_r$ is the dual  basis of $\g_1,\dots,\g_g$ then 
the coefficients $\U_j\bigl(s(\e_l)\bigr)$ in the 
description of $s_i$ are just the entries $c_{jl}$.

We have the following well known  property of the  
complex $B_\bullet(s, V)$. 

\begin{proposition}\label{P:Ei97} 
When the rank $r$ of $V$ is greater than or equal to the rank $e$ of 
$U$ the complex $B_\bullet(s,V)$ is exact if and only if 
the map $s$ is injective. When $r < e$ the complex  $B_\bullet(s, V)$ 
is exact if and only if $V=\im(s)$.  
\end{proposition}

\begin{proof} 
The proposition is immediate from Proposition \ref{P:koszulexact}, 
and  \cite{Ei97}, Theorem  A2.10.c. 
\end{proof} 

When $V=\im(s)$ we write $B_\bullet(s)$ instead of 
$B_\bullet(s,V)$. The exact complex $B_\bullet(s)$ is called the 
\emph{Buchsbaum-Rim} complex of the map $s$. 

Finally, we assign a multigrading on the components $B_i$ of 
the Buchsbaum-Rim complex $B_\bullet(s)$. The spaces 
$B_0=W$  and $B_1=U$ have the multigrading induced by the 
multidegrees of their basis  elements 
$\g_1,\dots, \g_g$ and $\e_1,\ldots, \e_e$, respectively.  
We set the multidegree of 
$\e_I=\e_{i_1}\wedge\dots\wedge\e_{i_p}$ 
in $\wedge^p U$ to be  $|\e_I|=\lcm(|\e_{i_1}|,\dots, |\e_{i_p}|)$,  
and thus obtain a multigrading 
on each of the vector spaces $\wedge^p U$. We also assign the 
multidegree $0$ to all elements of $D_{i-2}V^*$ and  
$\wedge^r V^*$. This way we get a multigrading of 
$B_i=D_{i-2}V^*\otimes\wedge^{i+r-1}U\otimes\wedge^r V^*$ 
for $i\ge 2$.

\section{Buchsbaum-Rim-Taylor theory}\label{S: BRTtheory}  

Let $s\: U\lra W$ be the map associated to the multigraded map 
$\phi: E\lra G$ as described in Section \ref{S: prelim}, and 
let $r=\rank \phi$.   
In this section we introduce the notion of a 
{\it BRT system} and show how a BRT system gives rise to 
a {\it BRT complex}\ :  
a finite free complex of multigraded $Q$-modules. 
This allows us to generalize to the case of 
multigraded modules the Taylor resolution \cite{Ta60}.     

Recall that $\Delta$ is the full simplex on the vertices 
$\{1,\dots, e\}$.  

\begin{definition} 
A family of vector spaces $\mathbb F=\{F_I\}_{I\in\Delta}$ is a 
{\it Buchsbaum-Rim-Taylor (BRT) system} for the map $s$ 
if the following three conditions are satisfied: 
\begin{enumerate} 
\item  $F_I=0$ whenever $|I|\le r$; 
\item  $F_I\subseteq D_{|I|-r-1}V^*$ whenever $|I|\ge r+1$; \ and   
\item  $\mathbb F$ is $s$-compatible: whenever  $|I|=p\ge r+2$, then 
\[ 
s_{p-r+1}(F_I\otimes \e_I\otimes\wedge^r V^*)\subseteq 
       \bigoplus_{\stackrel{\text{\scriptsize $|J|=p-1$}}{J\subset I}} 
          F_J\otimes\e_J\otimes\wedge^r V^*. 
\] 
\end{enumerate}
\end{definition} 

\begin{example} 
The main example of a BRT system is the \emph{full BRT system}  
$\mathbb F^{\text{full}}$ where  
\[
F_I= 
\begin{cases} 
D_{|I|-r-1}V^*  &\text{ if } |I|\ge r+1, \\ 
0               &\text{ otherwise}. 
\end{cases}
\]
The full BRT system is maximal, in the sense 
that it contains every other BRT system. 
In Section \ref{S: minimalres} we give 
another important example of a BRT system, the Scarf system. 
\end{example}

Having a BRT system $\mathbb F$ allows us to construct
a complex $R_\bullet(\mathbb F, \phi)$ of multigraded $Q$-modules  
as follows.

\begin{definition} 
We set $R_0=Q\otimes B_0=G$, and $R_1=Q\otimes B_1=E$. For $i\ge 2$ 
we define the multigraded $Q$-module 
\[ 
R_i= \bigoplus_{|I|=r+i-1} Q\otimes F_I\otimes 
                                   \e_I\otimes\wedge^r V^*.  
\]
We set $\phi_1=\phi$. For $i\ge 2$ 
we define the differentials $\phi_i\: R_i\lra R_{i-1}$ by 
homogenizing the restrictions of the maps $s_i$ to the 
free modules specified by the BRT system. More precisely,   
if $z\in F_I\otimes\e_I\otimes\wedge^r V^*$ and  $y\in Q$, and 
if $s_i^{I,J}$ is the component of 
$s_i$ that sends  $F_I\otimes \e_I\otimes \wedge^r V^*$ to 
$F_J\otimes\e_J\otimes\wedge^r V^*$, then the corresponding component 
$\phi_i^{I,J}$ of $\phi_i$ is 
\[ 
\phi_i^{I,J}(y\otimes z)=x^{|\e_I|-|\e_J|}y\otimes s_i^{I,J}(z).  
\] 
Since each map $\phi_i$ is obtained from the map $s_i$ by 
adjusting the multidegrees, and  
$B_\bullet(s)$ is a complex, it follows that 
\[ 
\begin{CD} 
R_\bullet(\mathbb F, \phi)= 0 \lra   
R_{e-r+1} @> \phi_{e-r+1} >>  R_{e-r} @>>>        \dots 
          @> \phi_2       >>  R_1     @> \phi >>  R_0 \lra 0  
\end{CD} 
\] 
is also a complex.  We call  $R_\bullet(\mathbb F, \phi)$ the 
{\it Buchsbaum-Rim-Taylor (BRT) complex} for the system $\mathbb F$ 
and the map $\phi$. 
\end{definition} 

\begin{definition}\label{D:genTaylor} 
We write $T_\bullet(\phi)$ for the complex 
$R_\bullet(\mathbb F^{\text{full}}, \phi)$,  
where $\mathbb F^{\text{full}}$ is the full BRT system, 
and call it the \emph{Taylor complex} of the map $\phi$. 
\end{definition} 

Note that if $J$ is a monomial ideal 
and $\phi$ is the minimal presentation map for the module $Q/J$, 
then $T_\bullet(\phi)$ is precisely the Taylor resolution 
\cite{Ta60} of $Q/J$.

We give a necessary and sufficient condition for the exactness 
of the complex $T_\bullet(\phi)$. For a multidegree $\a$ we say 
that the  coefficient matrix $C_{\a}$ is of \emph{maximal 
rank} if  $\rank C_\a=\min(r, \rank U_\a)$, that is, 
if the rank of $C_\a$ is the smaller of the rank of $\phi$ and 
the number of columns of $C_\a$.

\begin{theorem}\label{T: acyclicity}  
The complex $T_\bullet(\phi)$ is exact if and only if for 
every multidegree $\a$ the matrix $C_\a$ is of maximal rank. 
\end{theorem} 

\begin{proof} 
First we remark that a complex of multigraded $Q$-modules is exact 
if and only if it is exact in every multidegree $\a$.  

Next we notice that the multihomogeneous element 
$x^{\nu}\otimes v^{(\b)}\otimes\e_I\otimes v_{[r]}$ in $R_i$ is of 
multidegree $\a$ if  and only if the multidegree $\nu$ 
of the monomial $x^{\nu}$ added to the multidegree of 
$\e_I$  equals $\a$. In other words, $\e_I$ contributes to the 
component of multidegree $\a$ precisely when $|\e_i|\le \a$ for 
each $i\in I$. Therefore the component of $T_\bullet(\phi)$ 
of multidegree $\a$ is canonically isomorphic to the  
complex $B_\bullet(s_\a, V)$, where $V=\im(s)$, therefore 
by Proposition \ref{P:Ei97} is exact if 
and only if the matrix $C_\a$ is of maximal rank.  
\end{proof}

\begin{definition} 
We say that a multigraded morphism $\phi$ of rank $r$ is of 
\emph{uniform rank} if all $g\times r$ submatrices of 
its coefficient matrix $C$ have rank equal to $r$. 
\end{definition} 

The condition of Theorem \ref{T: acyclicity} is of course guaranteed 
whenever $\phi$ is of uniform rank. Thus we have 

\begin{corollary}\label{C:genericity} 
If $\phi$ is of uniform rank, then $T_\bullet(\phi)$ is exact. 
\qed
\end{corollary}

\begin{remark} 
The notion of uniform rank provides us with a precise description  
of what we mean when we say that a map $\phi$ is  
``sufficiently generic"  from the point of view of 
linear algebra.  Thus Corollary \ref{C:genericity} states that 
if $\phi$ is sufficiently generic from the point of view of 
linear algebra, then the exactness of the Taylor 
complex does not depend on the combinatorial structure of $\phi$ 
and the choices on the multidegrees of the generators of $E$ and $G$. 
\end{remark}

We conclude this section with an example.

\begin{example}\label{Ex: acyclicity}
Let $E=Q^4$ with standard basis $\e_1,\dots,\e_4$ of multidegrees 
\[
|\e_1|=(3,0), \quad |\e_2|=(2,1), \quad |\e_3|=(1,2), \quad   
|\e_4|=(0,3),  
\] 
and let $G=Q^2$ with standard basis of multidegrees 
$|\g_1|=(0,0)$, and $|\g_2|=(1,0)$. 

Let $\phi\: E \lra G$ be the multigraded homomorphism with 
standard matrix 
\[
\Phi= \left( 
\begin{array}{cccc} 
x^3 & x^2y & xy^2&y^3  \\ 
x^2 & 2xy &3y^2 &0 
\end{array} 
\right). 
\]  
Thus the coefficient matrix of $\phi$ is 
\[
C={\left (
\begin{array}{cccc}
1&1&1&1\\
1&2&3&0 
\end{array}\right)}, 
\] 
we have $r=\rank \phi =\rank s= 2$, and $V=\im s =W$.  
For the basis $\U_1, \U_2$ of $V^*=W^*$ we choose the dual 
of the standard basis $\g_1, \g_2$ of $W$.
Then the Taylor complex $T_\bullet(\phi)$ is  
\[ 
\begin{CD} 
0 \lra 
T_3 @>{\left( 
\begin{array}{cc} 
 y & 0  \\ 
-1 & -3 \\ 
1  & 2  \\ 
-x & -x 
\end{array} \right) } >>  
T_2 @>{\left( 
\begin{array}{cccc}
y^2  & -2y^3 & -3y^3 & 0      \\
-2xy & xy^2  & 0     & -3y^2  \\
x^2  & 0     & x^2y  & 2xy    \\
0    & x^3   & 2x^3  & x^2 
\end{array} \right) }>>       E \mapright{\Phi} G \lra 0;   
\end{CD} 
\] 
where the multidegrees of the generators of 
$T_2=Q\otimes\wedge^3 U\otimes \wedge^2 V^*\cong Q^4$ are  
\[ 
\begin{aligned}  
|1\otimes \e_{\{1,2,3\}}\otimes \U_{\{1,2\}}| &=(3,2),  \\ 
|1\otimes \e_{\{1,2,4\}}\otimes \U_{\{1,2\}}| &=(3,3),  \\ 
|1\otimes \e_{\{1,3,4\}}\otimes \U_{\{1,2\}}| &=(3,3),  \\  
|1\otimes \e_{\{2,3,4\}}\otimes \U_{\{1,2\}}| &=(2,3), 
\end{aligned}    
\]
and the multidegrees of the two generators of 
$T_3=Q\otimes V^*\otimes\wedge^4 U\otimes\wedge^2 V^* \cong Q^2$ are  
\[
\begin{aligned} 
|1\otimes \U_1 \otimes e_{\{1,2,3,4\}}\otimes \U_{\{1,2\}}| 
                            &= (3,3), \text{ and } \\ 
|1\otimes \U_2 \otimes e_{\{1,2,3,4\}}\otimes \U_{\{1,2\}}| 
                            &= (3,3). 
\end{aligned} 
\]  
Since $\phi$ is clearly of uniform rank (all $2\times 2$ minors of 
$C$ are non-zero), the complex   
$T_{\bullet}(\phi)$ is a resolution of $\coker(\phi)$. 
\end{example}

\section{The Scarf complex of a multigraded map}
\label{S: minimalres}

In this section we introduce and study the notion of a generic   
multigraded map and construct the minimal resolution of the  
cokernel of a minimal generic multigraded map. 

\begin{definition} 
Let $\phi\: E \lra G$ be a multigraded map. 
\begin{enumerate}
\item  The map $\phi$ is \emph{combinatorially generic} if 
       for every $1\le i<j\le e$ the support of $|\e_i|-|\e_j|$ 
       contains the supports of $|\e_i|$ and $|\e_j|$.  
\item  The map $\phi$ is \emph{generic} if it is combinatorially 
       generic and of uniform rank.  
\end{enumerate}
\end{definition} 

\begin{remark} 
(a) The notion of a combinatorially generic map is a translation to 
maps of the notion of a generic monomial ideal \cite{BaPeSt98}. When 
we deal with maps, we need to consider also the underlying linear 
algebra structure, hence the requirement that a generic map 
should be generic from both combinatorial and linear algebra 
point of view. 

(b) In \cite{MiStYa00} a different notion of generic monomial ideal 
is defined. It is straightforward to use it to define a different 
notion of generic multigraded map. However, since at this point 
we do not know whether our main result, Theorem \ref{T:Scarfres}, 
holds with this different definition of generic, we have elected 
not to pursue this line of investigation in this paper.  

(c) It is easy to see that via the deformation process of 
\cite{BaPeSt98} every map of uniform rank can be deformed to 
a generic map. We will see in this and the remaining sections that 
several important homological properties of monomial ideals hold 
also for maps of uniform rank. 
\end{remark}

Our next goal is to define the Scarf complex of 
a multigraded map $\phi$. 
This is a complex of multigraded free modules that 
is contained in the Taylor complex of $\phi$. It is 
a minimal complex when $\phi$ is minimal. In particular, it is 
contained in the minimal resolution of $M=\coker(\phi)$ when 
$\phi$ is minimal and of uniform rank. 
The main result of this paper, Theorem \ref{T:Scarfres},  
shows that the Scarf complex of $\phi$ is actually the minimal 
free resolution of $M$ when the map $\phi$ is minimal and generic.

We begin by translating to multigraded maps the notions of 
Scarf simplicial complex \cite{BaPeSt98}, and LCM-lattice 
\cite{GaPeWe99}. 

\begin{definition} 
Let $\phi\: E \lra G$ be a multigraded map, and let 
$\e_1,\dots,\e_e$ be a multihomogeneous basis of $E$.  
\begin{enumerate} 
\item  The \emph{Scarf simplicial complex} of $\phi$ 
       is the subcomplex 
       $\Delta_S=\Delta_S(\phi)$ of the full simplex $\Delta$ on 
       the vertices $\{ 1,\ldots, e\}$ defined as: 
\[ 
\Delta_S = \{ I\in \Delta\mid |\e_I|\ne |\e_J| 
                         \text{ for } J\ne I \}. 
\] 

\item  The \emph{LCM-lattice} $L_\phi$ of the map $\phi$ is the set  
       of all elements of $\mathbb N^n$ that can be obtained 
       as joins of some of the elements  
       $|\e_1|,\dots , |\e_e|$. Thus 
\[
L_\phi = \{ \a \mid  \a=|\e_I| \text{ for some face } I\in\Delta \}.  
\]
       is the set of those multidegrees that occur as multidegrees 
       of faces of $\Delta$. 
\end{enumerate} 
\end{definition}


Let $\phi\: E\lra G$ be a muiltigraded map of rank $r$, 
let $s\: U\lra W$ be defined as in Section \ref{S: prelim}, 
and let $V=\im(s)$. 
Before we define the Scarf system for $\phi$, we need to introduce 
some notation. 

We  partition  the LCM-lattice $L_\phi$ 
into two subsets: 
the subset $L_S$ containing the multidegrees of the faces of 
$\Delta_S$,  and its complement $L_S^0$. 
If $\a$ is a multidegree, we will 
denote by $I_\a$ the set of all indices $i$ for which 
$|\e_i|\preceq \a$.
Let $I(\a)$ be the intersection of all faces of $\Delta$ of 
degree $\a$, and set $I^\a = I_\a \smsm I(\a)$. Finally, 
recall from Section \ref{S: prelim} that $K_{I^\a}$ is the kernel 
of the canonical surjection $V^*\lra V^*_{I^\a}$.

\begin{definition}\label{D:Scarfsystem} 
The \emph{Scarf system} $\mathbb F_S$ for the map $\phi$ is 
the collection $\{ F_I\}$ of vector spaces defined by:  
\[
F_I= 
\begin{cases} 
0  &\quad\text{ if } |I|\le r; \\ 
D_{|I|-r-1}V^* 
   &\quad\text{ if $|I|\ge r+1$ and } I\in\Delta_S; \\  
D_{|I|-r-1}K_{I^\a} 
   &\quad\text{ if $|I|\ge r+1$ and $I=I_\a$ for some } \a\in L_S^0; \\ 
0  &\quad\text{ otherwise. }  
\end{cases} 
\]
\end{definition}

By Proposition \ref{P: ScarfBRT} the Scarf system is a 
BRT system, which allows for the following definition.

\begin{definition} 
We write $S_\bullet(\phi)$ for the BRT complex associated with 
the Scarf system $\mathbb F_S$ and call it the \emph{Scarf complex} 
of $\phi$.  
\end{definition}

The next theorem is the main result of this paper.

\begin{theorem}\label{T:Scarfres}
Let $\phi\: E \lra G$ be a minimal free multigraded presentation of 
a Noetherian multigraded $Q$-module $M$. If $\phi$ is generic, 
then the Scarf complex $S_\bullet(\phi)$ is a minimal free 
multigraded resolution of $M$ over $Q$. 
\end{theorem} 

We postpone the proof till the next section. We conclude this 
section with the proof that the Scarf system is a BRT system.

\begin{proposition}\label{P: ScarfBRT} 
The Scarf system is a BRT system. 
\end{proposition}

\begin{proof} 
We need to show that $\mathbb F_S$ is $s$-compatible. In other words 
if $I\in \Delta$ with $|I|=p\ge r+2$ we need to show that 
\[
s_{p-r+1}(F_I\otimes\e_I\otimes\wedge^r V^*)\subseteq 
\bigoplus_{\stackrel{\text{\scriptsize $|J|=p-1$}}{J\subset I}} 
          F_J\otimes\e_J\otimes\wedge^r V^* 
\] 
This is clear in all cases except when $I=I_\a$ for some 
$\a\in L^0_S$, so we assume this is the case. 
Thus $F_I=D_{p-r-1}K_{I^\a}$. Note that  
the component of 
$s_{p-r+1}$ in $F_J\otimes\e_J\otimes\wedge^r V^*$ is zero when 
$|\e_J|=\a$. Therefore, it will be enough to show that we have 
an inclusion 
\[
D_{p-r-2}K_{I^\a}\subseteq F_J 
\] 
whenever $J\subset I$ with $J=p-1$, and $|\e_J|=\b\ne \a$. 
This is clear if $J\in\Delta_S$, so we assume that $J\notin\Delta_S$. 
Then $J\subseteq I_\b\subsetneq I_\a=I$, hence $J=I_\b$. Thus it will 
be enough to show that $K_{I^\a}\subseteq K_{I^\b}$.   

Let $L$ be a face of $I_\b$ of  
degree $\b$. Let $\{i\}=I_\a\smsm I_\b$. 
Since $|\e_L|=\b$, we have that 
$|\e_L\wedge\e_i|=\a$. Thus  
$L\cup\{i\}$ contains $I(\a)$,   
hence $I(\a)\subseteq I(\b)\cup\{i\}$. 
Since $I_\b=I_\a \smsm \{i\}$, 
it follows that $I^\b \subseteq I^\a$, and therefore 
$K_{I^\a}\subseteq K_{I^\b}$, yielding the desired conclusion. 
\end{proof}

\begin{example}\label{Ex: Scarfres} 
Let $\phi$ be the multigraded map of Example \ref{Ex: acyclicity}.

For the LCM-lattice of $\phi$ and its parts $L_S$ and $L^0_S$ we have  
\[
\begin{aligned} 
L_\phi &= \{ \ (3,0), (2,1), (1,2), (0,3), (3,1), (3,2), (3,3), (2,2), 
               (2,3), (1,3) \ \};                                    \\ 
L_S    &= \{ \ (3,0), (2,1), (1,2), (0,3), (3,1), (2,2), (1,3) \ \}; \\ 
L^0_S  &= \{ \ (3,2), (2,3), (3,3) \ \}. 
\end{aligned} 
\] 

For the Scarf simplicial complex $\Delta_S$ we have  
\[
\Delta_S=\bigl\{ \{1\}, \{2\}, \{3\}, \{4\}, \{1,2\}, \{2,3\}, \{3,4\} 
         \bigr\} 
\] 

For $\a\in L^0_S$ the faces $I_\a$, $I(\a)$, and $I^\a$ are 
\[
\begin{aligned} 
I_{(3,2)} &= \{1,2,3\}, &  I(3,2) &= \{1,3\},  &  I^{(3,2)} &= \{2\}; \\ 
I_{(2,3)} &= \{2,3,4\}, &  I(2,3) &= \{2,4\},  &  I^{(2,3)} &= \{3\}; \\ 
I_{(3,3)} &= \{1,2,3,4\},\quad & I(3,3)&= \{1,4\},\quad & I^{(3,3)} &=\{2,3\}. 
\end{aligned} 
\]

For the spaces $K_{I^\a}$ with $\a\in L^0_S$ we obtain 
\[
\begin{aligned} 
K_{I^{(3,2)}}&=K_{\{2\}}   && =\Ker\bigl(V^* \lra V_{\{2\}}^*\bigr) &&=
              \Bbbk\cdot(2\upsilon_1 -\upsilon_2); \\ 
K_{I^{(2,3)}}&=K_{\{3\}}   && =\Ker\bigl(V^* \lra V_{\{3\}}^*\bigr) &&= 
              \Bbbk\cdot(3\upsilon_1 -\upsilon_3); \\ 
K_{I^{(3,3)}}&=K_{\{2,3\}} && =\Ker\bigl(V^* \lra V^*_{\{2,3\}}\bigr) &&= 0. 
\end{aligned} 
\]

For the spaces $F_I$ of the Scarf system for $\phi$ we obtain 
\[
\begin{aligned} 
{ }& F_{\{1,2,3\}}   &&= D_0 K_{I^{(3,2)}} = D_0 K_{\{2\}} &&= \Bbbk; \\ 
{ }& F_{\{2,3,4\}}   &&= D_0 K_{I^{(2,3)}} = D_0 K_{\{3\}} &&= \Bbbk; \\
{ }& F_{\{1,2,3,4\}} &&= D_1 K_{I^{(3,3)}} = K_{\{2,3\}}   &&= 0;     \\ 
{ }& F_I             &&= 0                             &&\text{otherwise}.  
\end{aligned} 
\]

Thus for the Scarf complex $S_\bullet(\phi)$ we obtain 
\[ 
\begin{CD} 
0 \lra 
S_2 @>{\left( 
\begin{array}{cccc} 
y^2 &0     \\ 
-2xy &-3y^2\\ 
x^2&2xy    \\
0&x^2 
\end{array} \right) }>>       
E @>\Phi >>  G @>>> 0    
\end{CD} 
\]
with $S_0=G\cong Q^2$, with $S_1=E \cong Q^4$, and with 
\[
S_2= \ \ Q\otimes \Bbbk\otimes \e_{\{1,2,3\}} \ \oplus \ 
         Q\otimes \Bbbk\otimes \e_{\{2,3,4\}} \ \ \cong Q^2. 
\]
Since the map $\phi$ is generic, the Scarf complex is the 
minimal free resolution of $M=\coker(\phi)$.  
\end{example}

\section{The proof of Theorem \ref{T:Scarfres}}\label{S:Proof}

Theorem \ref{T:Scarfres} is an immediate consequence of 
the following slightly stronger statement. 

\begin{theorem}\label{T:Scarfacyclicity} 
Let $\phi\: E \lra G$ be a minimal multigraded presentation of 
a Noetherian multigraded  $Q$-module $M$. 
If $\phi$ is combinatorially generic and 
the coefficient matrix $C_\a$ is of maximal rank for every 
multidegree $\a$, then the Scarf complex $S_\bullet(\phi)$ 
is the minimal free resolution of $M$.
\end{theorem}  

\begin{proof} 
It is clear from the construction that the Scarf complex is 
minimal, hence it suffices to show that it is exact. Also, 
by construction $S_\bullet=S_\bullet(\phi)$ is a subcomplex of the 
Taylor complex $T_\bullet=T_\bullet(\phi)$, and $T_\bullet$ 
is exact by Theorem \ref{T: acyclicity}. Therefore it is 
enough to show that $X_\bullet=T_\bullet/S_\bullet$ is an 
exact complex. 
We will do this by showing that there is a filtration of $X _\bullet$
whose $i$th quotient $X^i_\bullet$ is the direct sum of exact
complexes.

Note that 
the component of $X_\bullet$ in homological degree $m\ge 2$ is   
\[ 
X_m = 
\bigoplus_{\stackrel{\text{\scriptsize $I\notin\Delta_S$}}{|I|=m+r-1}} 
Q\otimes  H_I \otimes\e_I\otimes\wedge^r V^*  
\]
where   
\[ 
H_I= 
\begin{cases} 
D_{m-2} V^*/ D_{m-2} K_{I^\a}  
     &\quad\text{ if $I=I_\a$ for some } \a\in L_S^0; \\ 
D_{m-2} V^*  
     &\quad\text{ otherwise};  
\end{cases} 
\] 
and $X_m=0$ for $m\le 1$.

We partition $L_S^0$ as follows.  
Let $L_S^1$ be the set of minimal multidegrees in $L_S^0$ 
(with respect to the partial order $\prec$).  
Once the sets $L_S^1,\dots, L_S^i$ 
have been defined, we define $L_S^{i+1}$ to be the set of minimal 
elements in $L_S^0 \smsm (L_S^1\cup\ldots\cup L_S^i)$. 
We also define $UL_S^i=L_S^1\cup\ldots\cup L_S^i$. 
Now define $X_m^i\subseteq X_m$ by 
\[
X_m^i = 
\bigoplus_{\stackrel
          {\text{\scriptsize $I(\a)\subseteq I\subseteq I_\a$}}
          {|I|=m+r-1; \ \a\in UL_S^i} } 
Q\otimes  H_I \otimes\e_I\otimes\wedge^r V^*.   
\]
It is straightforward that $X^i_\bullet$ is a subcomplex of 
$X_\bullet$ for each $i$, and that the quotient complex 
$Y^i_\bullet= X^i_\bullet / X^{i-1}_\bullet$ has as its component 
in homological degree $m\ge 2$ the module 
\[
Y^i_m=
\bigoplus_{\stackrel
          {\text{\scriptsize $I(\a)\subseteq I\subseteq I_\a$}}
          {|I|=m+r-1; \ \a\in L_S^i} } 
Q\otimes  H_I \otimes\e_I\otimes\wedge^r V^*,   
\]   
while $Y^i_m=0$ for $m=0,1$. Thus to show that $S_\bullet$ is exact, 
it is enough to show that the complex $Y^i_\bullet$ is split 
exact for each $i\ge 1$.  Note however that the complex 
$Y^i_\bullet$ decomposes into the direct sum of subcomplexes 
\[
Y^i_\bullet = \bigoplus_{\a \in L_S^i} Y^i(\a), 
\]
where the complex $Y^i(\a)$ has as its component in homological 
degree $m\ge 2$ the module  
\[
Y^i(\a)_m= 
\bigoplus_{\stackrel
          {\text{\scriptsize $I(\a)\subseteq I\subseteq I_\a$}}
          {|I|=m+r-1}} 
Q\otimes  H_I \otimes\e_I\otimes\wedge^r V^*,   
\]
and $Y^i(\a)_m=0$ for $m=0,1$. Thus it suffices to show that 
each complex $Y^i(\a)$ is split exact. For the rest of this 
proof we fix $i$ and $\a\in L_S^i$, and we write $Z_\bullet$ 
for the complex $Y^i(\a)$. 

Since $\phi$ is combinatorially generic, if $\a\in L_S^0$, and 
$I$ and $J$ are  such that $|\e_I|=|\e_J|=\a$, then
$|\e_{I\cap J}|=\a$ as well. Therefore  
there exists a  unique minimal face $I(\a)$ of degree $\a$.   
Furthermore, a face 
$I$ has degree $\a$ if and only if we have 
$I(\a)\subseteq I\subseteq I_\a$, while 
the containment $I(\a)\subset I_\a$ is strict. 
So  $i\in I_\a \smsm I(\a)$ if and only if  
$|\e_{I_\a \smsm i}|=|\e_{I_\a}|=\a$, and it follows that
$I^\a=I_\a \smsm I(\a)=I_{\a'}$ where $\a'=\a - (1,\dots,1)$. 
Thus $K_{I^\a}=K_{I_{\a'}}=K_{\a'}$, and the short exact sequence 
of vector spaces 
\[
0 \lra K_{\a'} \lra V^* \lra  V^*_{\a'} \lra 0 
\]
induces for $m\ge 2$ a canonical filtration 
\[
0=T^m_{-1}\subseteq T^m_0\subseteq \dots 
                  \subseteq T^m_{m-2}=D_{m-2}V^* 
\]
on $D_{m-2}V^*$ whose $i$th quotient $T^m_i/T^m_{i-1}$ is canonically 
isomorphic for $i\ge 0$ to $D_{m-2-i} K_{\a'}\otimes D_i V^*_{\a'}$. 
This way we obtain a filtration 
\[
0=P^{-1}_m(I)\subseteq 
P^0_m(I)\subseteq P_m^1(I)\subseteq\dots\subseteq P_m^k(I)\subseteq\dots 
\]
on each $H_I$, where for $i\ge 0$ we set   
\[
P^i_m(I)= 
\begin{cases} 
T^m_i  &\quad\text{if } I\ne I_\a; \\ 
T^m_i/T^m_0 &\quad\text{otherwise}. 
\end{cases} 
\]
When $i\ge 0$ it is clear that for the $i$th quotient of this filtration 
we have canonically 
\[
P^i_m(I)/P^{i-1}_m(I)\quad\cong\quad 
\begin{cases} 
0              &\quad\text{if $I=I_\a$ and $i=0$}; \\ 
D_{m-2-i}K_{\a'}\otimes D_i V_{\a'}^* 
               &\quad\text{otherwise.} 
\end{cases}  
\] 
Let $t=|I_\a|$. Set $Z^{-1}_m=0$, and for $k\ge 0$ 
define $Z^k_m\subseteq Z_m$ as 
\[
Z^k_m= 
\bigoplus_{\stackrel
          {\text{\scriptsize $I(\a)\subseteq I\subseteq I_\a$}}
          {|I|=m+r-1}} 
Q\otimes  P_m^{m+r+k-t}(I) \otimes\e_I\otimes\wedge^r V^*,   
\]
It is straightforward from these definitions that 
$Z^k_\bullet$ is a subcomplex of $Z_\bullet$ for each $k\ge -1$.  
Note that if $k\ge 0$ and $I=I_\a$ then  $m+r+k-t=k+1\ge 1$. 
Therefore for each $k\ge 0$ the quotient complex 
$\overline{Z^k_\bullet}= Z^k_\bullet/ Z^{k-1}_\bullet$ has as its 
component in homological degree $m\ge 2$ the module 
\[
\overline{Z^k_m} = 
\bigoplus_{\stackrel
          {\text{\scriptsize $I(\a)\subseteq I\subseteq I_\a$}}
          {|I|=m+r-1}} 
Q\otimes D_{t-r-k-2} K_{\a'}\otimes D_{m+r+k-t} V^*_{\a'}  
                    \otimes\e_I\otimes\wedge^r V^*,   
\]
and is $0$ in homological degrees $m=0,1$. Thus to complete the proof 
of the theorem it suffices to show that $\overline{Z^k_\bullet}$ 
is split exact.  

Let $q=|I(\a)|$. 
We examine the differential of $\overline{Z^k_\bullet}$.

When $k>t-q-1$ 
it is clear that we have a canonical isomorphism of complexes
\[
\overline{Z^k_\bullet}\cong 
Q\otimes D_{t-r-k-2} K_{\a'} 
 \otimes A^{q+k+1-t, \ 0}_\bullet(s_{\a'})[-q+r-1]   
 \otimes\e_{I(\a)}\otimes\wedge^r V^*.    
\] 
By Proposition \ref{P:koszulexact}  the complex 
$A^{q+k+1-t, \ 0}_\bullet(s_{\a'})$ is split exact, therefore
$\overline{Z^k_\bullet}$ is split exact as well.

Similarly when $k\le t-q-1$
we have a canonical isomorphism of complexes 
\[
\overline{Z^k_\bullet}\cong 
Q\otimes D_{t-r-k-2} K_{\a'} 
 \otimes A^{0,\ t-q-k-1}_\bullet(s_{\a'})[-p]   
 \otimes\e_{I(\a)}\otimes\wedge^r V^*,   
\] 
where the shift $p$ is computed as   
$p=\max\bigl(q-r+1,\ 2\bigr)$. 

Next, recall that by assumption the rank $r_{\a'}$ of $s_{\a'}$ 
is equal to $\min(r, |I_{\a'}|)$. 
If $|I_{\a'}|=r_{\a'}\le r$ then $t-q-k-1=|I_{\a'}|-k-1 < r_{\a'}$\ ; 
thus 
$A^{0,\ t-q-k-1}_\bullet(s_{\a'})$ (hence also $\overline{Z^k_\bullet}$) 
is  split exact by Corollary \ref{T:splitexact}.

Finally, assume $|I_{\a'}| > r$. Then $\rank s_{\a'}=r$, hence  
$K_{\a'}=0$. If $t-r-k-2\ne 0$ then 
$D_{t-r-k-2} K_{\a'}=0$, therefore $\overline{Z^k_\bullet}=0$ 
is split exact. 
If $t-r-k-2 =0$ then    
\[ 
t-q-k-1= r + 2 -q -1 = r+1-|I(\a)|< r = r_{\a'}, 
\]  
where the last inequality follows by Lemma \ref{L:sizeofI}. Therefore
$A^{0,\ t-q-k-1}_\bullet(s_{\a'})$ (hence also 
$\overline{Z^k_\bullet}$) is split exact by Corollary 
\ref{T:splitexact} 
\end{proof} 

\begin{lemma}\label{L:sizeofI} 
With the assumptions of  Theorem \ref{T:Scarfacyclicity}, if 
$|I_{\a'}| > r$ then $|I(\a)|\ge 2$.
\end{lemma}

\begin{proof} 
Suppose that $|I(\a)|=1$. Thus $I(\a)=\{l\}$ and 
$|e_l|=\a$. Since by asumption $C_{\a'}$ is 
of maximal rank, we can choose $I'\subset I_{\a'}$
such that $|I'|=r$ and $\rank s_{I'}=r$. Let $I=I'\cup\{l\}$.  Then  
$\phi_2\bigl(1\otimes \e_I\otimes\U_{[r]}\bigr)$ is a minimal 
free generator of $E$, and is a syzygy of $\phi$. This 
contradicts the assumption that $\phi$ is a minimal 
presentation of the multigraded module $M$.
\end{proof}

\section{The LCM-Lattice}\label{S: LCM}

The importance of the LCM-lattice in determining the minimal 
resolutions of monomial ideals was exhibited in \cite{GaPeWe99}. 
In the setting of multigraded maps one has to take into account 
also the underlying linear algebra structure. Our goal is to show 
that for maps of uniform rank, if the linear algebra 
structure is essentially the same then the structure of the minimal 
resolution is determined, subject to a certain compatibility condition,  
by the isomorphism class of the LCM-lattice.

Let $\phi\: E \lra G$ 
be a multigraded map of finite free $Q$-modules, 
let  $s\: U\lra W$ be the associated map of vector spaces (see 
Section \ref{S: prelim}), let $V$ be the image of $s$, 
and let $L_{\phi}$ be the LCM-lattice. Similarly, 
let $Q'$ be another polynomial ring, let $\phi'\: E'\lra G'$ 
be a multigraded map of finite free $Q'$-modules, and consider 
the corresponding objects $s'$, $V'$, and $L_{\phi'}$ . 
Let $e$ be the rank of $E$, and let $e'$ be the rank of $E'$. 

\begin{definition} 
(a)  The maps $\phi$ and $\phi'$ are called \emph{quasi-equivalent} if 
there exists a choice of homogeneous bases for $E$ and $E'$, 
and a choice of bases for $V$ and $V'$ such that the matrices 
of the induced by $s$ and $s'$ maps $U\lra V$ and $U'\lra V'$ 
are the same. 

(b) A choice of bases of $E, E', V$, and $V'$ as in (a) is called a 
\emph{QE-structure} for the pair $(\phi,\phi')$. 

(c) Assume $\phi$ and $\phi'$ are quasi-equivalent. A function of sets  
$f\: L_\phi\lra L_{\phi'}$ is called \emph{QE-compatible} if there 
exist multihomogeneous bases $\e_1,\dots,\e_e$ of $E$ and 
$\e'_1,\dots,\e'_e$  of $E'$ that are part of a QE-structure 
and satisfy $f(|\e_i|)=|\e'_i|$ for all $i$. 
\end{definition}

To state our theorem, we need to introduce a slight generalization 
of  the relabeling procedure of \cite{GaPeWe99}. 
Let $T_\bullet$ be an extension of $\phi$ to a  
finite free multigraded complex of the form 
\[
0\lra T_p \mapright{\phi_p} T_{p-1} \lra \dots \lra 
      T_2 \mapright{\phi_2} T_1 \mapright{\phi} T_0\lra  0   
\]
where $T_0=G$ and $T_1=E$; such that the free modules 
$T_m$ for $m\ge 1$ have generators with multidegrees in $L_{\phi}$.  
Let $f\: L_{\phi}\lra L_{\phi'}$ be a QE-compatible map  
which preserves joins that appear as 
multidegrees of minimal generators of the free modules $T_m$ for 
$m\ge 1$. 
                                                                               
Using the map $f$ we relabel  $T_\bullet$ in the spirit of   
\cite{GaPeWe99}, Construction 3.2.  First, we replace $T_0$ with 
$G'$, and $\phi$ with $\phi'$. Next,  if a free copy of  
$Q$ in $T_m$ for some $m\ge 1$ has multidegree $\a$  then after relabeling 
we get a free copy of $Q'$ with  multidegree $f(\a)$. Finally, 
we relabel the maps $\phi_i$ for $i\ge 2$ to $\phi'_i=f(\phi_i)$  
by homogenizing  the images of $\phi_i$. We write $f(T_\bullet)$ 
for the resulting complex.

We are now ready to state 
the generalization of \cite{GaPeWe99}, Theorem 3.3. 

\begin{theorem}\label{T: lcmres} 
Let $\phi$ and $\phi'$ be two quasi-equivalent multigraded 
maps  of uniform rank $r$. 
Let $f\:  L_{\phi}\lra L_{\phi'}$ be a QE-compatible map which 
preserves joins of any $s$ atoms, where $s\ge r+1$. 
Let $F_\bullet (\phi)$ be the minimal resolution of the cokernel of 
the map $\phi$. 
Then $f(F_\bullet (\phi))$ is a free resolution of the 
cokernel of $\phi'$. 
\end{theorem} 

\begin{proof}
By Corollary \ref{C:genericity}, the Taylor complexes  
$T_\bullet(\phi)$ and $T_\bullet(\phi')$ 
resolve the \hyphenation{co-kernels} 
cokernels of $\phi$ and $\phi'$.  Since  
$f(T_\bullet(\phi))=T_\bullet(\phi')$, the zeroth homology of 
$f(T_\bullet(\phi))$ is equal to the cokernel of $\phi'$. Because  
$T_\bullet(\phi)=F_\bullet (\phi)\oplus P_\bullet$ where $P_\bullet$ 
is a direct sum of split exact complexes
of the form $0 \lra  Q\lra Q \lra 0$, it follows that 
$f(T_\bullet(\phi))=f(F_\bullet (\phi)) \oplus f(P_\bullet)$ where
$f(P_\bullet)$ is the direct sum of split exact complexes
of the form $0 \lra  Q'\lra Q' \lra 0$. 
\end{proof} 

The remarks made in \cite{GaPeWe99}, Example 3.4, hold for the obvious
generalizations for multigraded maps. 
We finish this section with an example where we apply 
Theorem \ref{T: lcmres}. 

\begin{example}\label{Ex: lcm} 
Let $Q'=\Bbbk[u,v,w]$, and 
let $\phi'\: (Q')^4 \lra (Q')^2$ be the multigraded homomorphism 
with standard matrix: 
\[
 \Phi'=
 \left( 
\begin{array}{cccc} 
u^2v & uvw & u^2w &uw^2  \\ 
uv & 2vw &3uw &0 
\end{array} \right). 
\] 
The coefficient matrix of $\phi'$ equals the coefficient matrix of
$\phi$ of Examples \ref{Ex: acyclicity} 
and  \ref{Ex: Scarfres}, hence the maps $\phi$ and $\phi'$ are 
quasi-equivalent. Note that $\phi'$ is not a generic
map.  Define $f\: L_{\phi}\lra L_{\phi'}$ by 
\[
\begin{aligned} 
f(3,0) &=(2,1,0),\qquad & f(2,0)&=(1,1,1),\qquad & f(1,2)&=(2,0,1), \\ 
f(0,3) &=(1,0,2),       & f(3,2)&=(2,1,1),       & f(3,3)&=(2,1,2), \\ 
       &                & f(2,3)&=(2,1,2).       &
\end{aligned} 
\]
Then $f$ satisfies the conditions of Theorem \ref{T: lcmres}, 
but is not an isomorphism of LCM-lattices. Thus 
a free resolution of  the cokernel of $\phi'$  
can be obtained by applying  
Theorem \ref{T: lcmres} to the minimal resolution of 
Example \ref{Ex: Scarfres}: 
\[ 
\begin{CD} 
0 \lra 
 (Q')^2 @>{\left( 
\begin{array}{cccc} 
w &0      \\ 
-2u &-3uw \\ 
v&2v      \\ 
0&uv
\end{array} \right) }>>       
 (Q')^4 @>\Psi >>  (Q')^2 @>>> M' @>>>0.  
\end{CD} 
\] 
\end{example}

\end{document}